\documentclass[12pt,oneside]{amsart}

\newtheorem{thm}{Theorem}[section] 
\newtheorem{cor}[thm]{Corollary}

\newtheorem{lem}[thm]{Lemma} 

\newtheorem{prop}[thm]{Proposition}

\theoremstyle{definition}
\theoremstyle{remark}
\theoremstyle{proof}

\numberwithin{equation}{section}

\newcommand{\norm}[1]{\left\Vert#1\right\Vert}

\newcommand{\set}[1]{\left\{#1\right\}}
\newcommand{\brac}[1]{\left(#1\right)}
\newcommand{\scalar}[1]{\left \langle #1 \right \rangle}

\newcommand{\Real}{\mathbb{R}}

\newcommand{\eps}{\varepsilon}

\newcommand{\K}{\mathcal{K}}

\def \eps {\varepsilon}

\begin{document}

\title{A remark on two duality relations}

\author{Emanuel Milman \\
\fontfamily{cmr} \fontseries{l} \fontshape{sc} \fontsize{12}{0}
\selectfont \medskip \today}
\thanks{Supported in part by BSF and ISF}


\maketitle

\begin{abstract}
We remark that an easy combination of two known results yields a
positive answer, up to $log(n)$ terms, to a duality conjecture
that goes back to Pietsch. In particular, we show that for any two
symmetric convex bodies $K,T$ in $\Real^n$, denoting by $N(K,T)$
the minimal number of translates of $T$ needed to cover $K$, one
has:
\[
\!\!\!\!\!\!\!\!\!\!\!\!\! N(K,T) \leq N(T^\circ,(C \log(n))^{-1}
K^\circ)^{C \log(n) \log\log(n)},
\]
where $K^\circ,T^\circ$ are the polar bodies to $K,T$,
respectively, and $C\geq 1$ is a universal constant. As a
corollary, we observe a new duality result (up to $\log(n)$
terms) for Talagrand's $\gamma_p$ functionals.
\end{abstract}


\section{Introduction}

Let $K$ and $T$ denote two convex bodies in $\Real^n$ (i.e.
convex compact sets with non-empty interior). Throughout this
note we assume that all bodies in question are centrally
symmetric w.r.t. to the origin (e.g. $K=-K$). For a convex body
$L$, we denote by $L^\circ$ its polar body, defined as $L^\circ =
\set{ x \in \Real^n ; \scalar{x,y} \leq 1 \; \forall y \in L }$.
The \emph{covering number of $K$ by $T$}, denoted $N(K,T)$, is
defined as the minimal number of translates of $T$ needed to
cover $K$, i.e.:
\[
N(K,T) = \min\set{ N\; ; \; \exists x_1,\ldots,x_N \in \Real^n
\;,\; K \subset \bigcup_{1\leq i \leq N} (x_i + T) }.
\]
In this note, we address the following conjecture of Pietsch
(\cite[p. 38]{Pietsch-Book}) from 1972, originally formulated in
operator-theoretic notations:

\medskip
\noindent\textbf{Duality Conjecture for Covering Numbers.} Do
there exist numerical constants $a, b \geq 1$ such that for any
dimension $n$ and for any two symmetric convex bodies $K,T$ in
$\Real^n$ one has:
\begin{equation} \label{eq:duality-covering}
b^{-1} \log N(T^\circ,a K^\circ) \leq \log N(K,T) \leq b \log
N(T^\circ,a^{-1} K^\circ) \;\; ?
\end{equation}

This problem may be equivalently formulated using the notion of
\emph{entropy numbers}. For a real number $k \geq 0$, denote the
$k$'th entropy number of $K$ w.r.t. $T$ as:
\[
e_k(K,T) = \inf \{ \eps>0 ; N(K,\eps T) \leq 2^k \}.
\]
Then the duality conjecture may be equivalently formulated with
(\ref{eq:duality-covering}) replaced by:
\begin{equation} \label{eq:duality-entropy}
a^{-1} e_{bk}(T^\circ,K^\circ) \leq e_k(K,T) \leq a e_{b^{-1}
k}(T^\circ,K^\circ)
\end{equation}
for all $k\geq 0$ (and there is no loss in generality if we
assume that $k$ is an integer).

As already mentioned, the duality conjecture originated from
operator theory, where entropy numbers are used to quantify the
compactness of an operator $u : X \rightarrow Y$ between two
Banach spaces. Leaving the finite dimensional setting for a brief
moment, if $K = u(B(X))$ and $T = B(Y)$, where $B(Z)$ denotes the
unit-ball of a Banach space $Z$, then it is easy to see that
$e_k(K,T) \rightarrow 0$ as $k\rightarrow \infty$ iff the operator
$u$ is compact. Since $u$ is compact iff its dual $u^*:Y^*
\rightarrow X^*$ is too, and since $u^*(B(Y^*)) = u^*(T^\circ)$
and $B(X^*) = u^*(K^\circ)$, it follows that $e_k(K,T) \rightarrow
0$ iff $e_k(T^\circ,K^\circ) \rightarrow 0$. Hence, it is natural
to conjecture that the rate of convergence to 0 is asymptotically
similar in both cases. A strong interpretation of this similarity
is given by (\ref{eq:duality-entropy}). We will mention other
weaker interpretations below.

\medskip

Although the general problem is still not completely settled,
there has been substantial progress in recent years, and the
answer is known to be positive for a wide class of bodies. We
begin by describing some results in this direction. We comment
here that when the result imposes the same restrictions on $K$
and $T$, it is obviously enough to specify only one side of the
inequalities in (\ref{eq:duality-covering}) or
(\ref{eq:duality-entropy}). When both $K$ and $T$ are ellipsoids,
it is easy to see that in fact $N(K,T) = N(T^\circ,K^\circ)$.
Other special cases were settled in
\cite{Schutt-Entropy-numbers},\cite{Carl-Entropy-numbers},\cite{GKS-Entropy-numbers},
\cite{KMT-Entropy-numbers},\cite{PajorTomczak-s-numbers}. In
\cite{Konig-Milman}, it was shown that:
\[
C^{-n} N(T^\circ,K^\circ) \leq N(K,T) \leq C^{n}
N(T^\circ,K^\circ),
\]
for some universal constant $C>1$. This implies that the tail
behaviour of the entropy numbers satisfies the duality problem,
i.e. $e_{\lambda k}(K,T) \leq 2 e_k(T^\circ,K^\circ)$ for some
universal constant $\lambda>0$ and all $k \geq n$. This was
subsequently generalized in \cite{Pisier-Regular-M-Position}.

Another variant of the problem, is to consider not the individual
entropy numbers, but rather the entire sequences $\set{e_k(K,T)}$
and $\set{e_k(T^\circ,K^\circ)}$. Then one may ask whether:
\begin{equation} \label{eq:duality-sequence}
C^{-1} \norm{\set{e_k(T^\circ,K^\circ)}} \leq
\norm{\set{e_k(K,T)}} \leq C \norm{\set{e_k(T^\circ,K^\circ)}}
\end{equation}
for some universal constant $C>1$ and any symmetric (i.e.
invariant to permutations) norm $\norm{\cdot}$. When one of the
bodies is an ellipsoid, this was positively settled in
\cite{Tomczak-Sequential-Duality-For-Ball}. Later, in \cite{BPST},
this was extended to the case when one of the bodies is uniformly
convex or more generally K-convex (see \cite{BPST} and
\cite{Pisier-Book} for definitions), in which case the constant
$C$ in (\ref{eq:duality-sequence}) depends only on the
K-convexity constant. The technique developed in \cite{BPST}
played a crucial role in some of the subsequent results on this
problem, and one particular remark will play an essential role in
this note.

Returning to the duality problem of individual entropy numbers,
it was shown in \cite{Milman-Szarek-Geometric-Lemma} that there
exist universal constants $a,b\geq 1$ such that when $T=D$ is an
ellipsoid:
\[
e_{bk}(D^\circ,K^\circ) \leq a (1+\log k)^3 e_k(K,D),
\]
for all $k \geq 0$. In addition, the authors of
\cite{Milman-Szarek-Geometric-Lemma} observed a connection
between (one side of) the duality conjecture with $T=D$ and a
certain geometric lemma. Later, the case when one of the bodies
is an ellipsoid was completely settled in
\cite{AMS-Duality-For-Ball}, by showing that:
\[
b^{-1} \log N(D^\circ,a K^\circ) \leq \log N(K,D) \leq b \log
N(D^\circ,a^{-1} K^\circ).
\]
The main new tool developed in \cite{AMS-Duality-For-Ball} was the
so called ``Reduction Lemma", which roughly reduces the problem
(\ref{eq:duality-covering}) for all $K,T$ to the case $K \subset
4 T$. This will be the second important tool in this note.

Finally, in \cite{AMST-Duality-For-K-Convex}, the Reduction Lemma
was combined with the techniques developed in \cite{BPST}, to
transfer the results obtained there for the sequence of entropy
numbers, to the individual ones. Thus, when one of the bodies $K$
or $T$ is K-convex, (\ref{eq:duality-covering}) was shown to hold
with the constants $a,b$ depending solely on the K-convexity
constant. The key ideological step in
\cite{AMST-Duality-For-K-Convex} was to separate the question of
``complexity" from the question of duality, by explicitly
introducing a new notion of \emph{convexified packing number},
which was implicitly used in \cite{BPST}. We will later refer to
this new notion as well.

\medskip

Our first new observation in this note is in fact an immediate
consequence of Theorem 6 in \cite{BPST} and the Reduction Lemma in
\cite{AMS-Duality-For-Ball}. It settles the duality problem
(\ref{eq:duality-covering}) (and (\ref{eq:duality-entropy})) up
to $\log(n)$ terms, and in fact strengthens and generalizes all
previously known results into a single statement. Because of the
symmetry between $K$ and $T$ (as explained below), we formulate
this as a one sided inequality:

\begin{thm} \label{thm:main}
Let $K,T$ be two symmetric convex bodies in $\Real^n$. Then:
\begin{equation} \label{eq:final}
\log N(K,T) \leq V \log(V) \log N(T^\circ,V^{-1} K^\circ),
\end{equation}
where $V = \min( V(K), V(T) )$ and $V(L)$ is defined as:
\begin{equation} \label{eq:defn-V}
V(L) := \inf \set{ \log(C d_{BM}(L,B)) f(\K(X_B)) ; B \text{ is a
convex body in } \Real^n },
\end{equation}
where $C>0$ is a universal constant, and $f$ is a function
depending solely on $\K(X_B)$, the K-convexity constant of the
Banach space $X_B$ whose unit ball is $B$.
\end{thm}

Recall that the Banach-Mazur distance $d_{BM}(L,B)$ of two
symmetric convex bodies $L,B$ is defined as:
\[
d_{BM}(L,B) := \inf \set{ \gamma \geq 1; B \subset T(L) \subset
\gamma B},
\]
where the infimum runs over all linear transformations $T$. Since
$V(L) = V(L^\circ)$ because $d_{BM}(L,B) =
d_{BM}(L^\circ,B^\circ)$ and $\K(X_{B^\circ}) = \K(X_B)$,
applying the Theorem to $K' = T^\circ$ and $T' = K^\circ$ gives
the opposite inequality (with the same $V$):
\begin{equation} \label{eq:other-side}
(V \log(V))^{-1} \log N(T^\circ,V K^\circ) \leq \log N(K,T).
\end{equation}
In addition, since by John's Theorem, the Banach-Mazur distance
of any symmetric convex body in $\Real^n$ from the Euclidean ball
$D$ is at most $\sqrt{n}$, and since $\K(D)=1$, we immediately
have:

\begin{cor} \label{cor:main}
With the same notations as in Theorem \ref{thm:main}:
\[
\log N(K,T) \leq C \log(1+n) \log \log(2+n) \log N(T^\circ,(C
\log(1+n))^{-1} K^\circ),
\]
where $C>0$ is a universal constant.
\end{cor}

This should be compared with the previously known best estimate
(to the best of our knowledge) for general symmetric convex bodies
$K,T$:
\[
\log N(K,T) \leq C \log N(T^\circ, (C n)^{-1/2} K^\circ),
\]
which is derived by comparing $K$ with its John ellipsoid and
using the duality result of \cite{AMS-Duality-For-Ball} for
ellipsoids.

\medskip

Although there has been much progress in recent years towards a
positive answer to the duality conjecture, it is still not clear
that a positive answer should hold in full generality. In view of
the Corollary \ref{cor:main}, and Pisier's well known estimate
$\K(X_B) \leq C \log(1+n)$ for any symmetric convex body $B$ in
$\Real^n$, we conjecture a weaker form of the duality problem:

\medskip
\noindent\textbf{Weak Duality Conjecture for Covering Numbers.}
Does there exist a numerical constant $C \geq 1$ such that for any
dimension $n$ and for any two symmetric convex bodies $K,T$ in
$\Real^n$ one has:
\[
\log N(K,T) \leq V \log N(T^\circ,V^{-1} K^\circ),
\]
where $V = C \min( \K(X_K), \K(X_T) )$ ?

\medskip

We present the proof of Theorem \ref{thm:main} and several other
connections to previously mentioned notions in Section
\ref{sec:1}. In Section \ref{sec:2}, we give an application of
Corollary \ref{cor:main} for Talagrand's celebrated $\gamma_p$
functionals, which was in fact our motivation for seeking a
result in the spirit of Corollary \ref{cor:main}. Recall that for
a metric space $(M,d)$ and $p>0$, $\gamma_p(M,d)$ is defined as:
\[
\gamma_p(M,d) := \inf \sup_{x \in M} \sum_{j \geq 0} 2^{j/p}
d(x,M_j)
\]
where the infimum runs over all \emph{admissible} sets
$\set{M_j}$, meaning that $M_j \subset M$ and $|M_j| = 2^{2^j}$
(we refer to \cite[Theorem 1.3.5]{Talagrand-Book} for the
connection to equivalent definitions). For two symmetric convex
bodies $K,T$, let us denote $\gamma_p(K,T) := \gamma_p(K,d_T)$,
where $d_T$ is the metric corresponding to the norm induced by
$T$. The $\gamma_2(\cdot,D)$ functional, when $D$ is an
ellipsoid, was introduced to study the boundedness of Gaussian
processes (see \cite{Talagrand-Book} for an historical account on
this topic). It was shown by Talagrand in his celebrated
``Majorizing Measures Theorem", that in fact $\gamma_2(K,D)$ and
$E \sup_{x\in K} \scalar{x,G}$, where $G$ is a Gaussian r.v. (with
covariance corresponding to $D$ in an appropriate manner), are
equivalent to within universal constants. This was later extended
to various other classes of stochastic processes, where the
naturally arising metric $d$ is not the $l_2$ norm (again we
refer to \cite{Talagrand-Book} for an account).

Our second observation in this note is the following duality
relation for the $\gamma_p$ functionals:

\begin{thm} \label{thm:gamma}
Let $K,T$ be two symmetric convex bodies in $\Real^n$. Then for
any $p>0$:
\[
\gamma_p(K,T) \leq C_p \log(1+n)^{2+1/p} \log \log(2+n)^{1/p}
\gamma_p(T^\circ,K^\circ),
\]
where $C_p > 0$ depends solely on $p$.
\end{thm}

Although we strongly feel that this is unlikely, one could
conjecture that the $\log(n)$ terms are not required in the last
Theorem (at least for some values of $p$). In that case, as will
be evident from the proof, we mention that such a conjecture is
independent of the duality conjecture for covering numbers, in
the sense that neither one implies the other.

\medskip

\noindent \textbf{Acknowledgments.} I would like to sincerely
thank my supervisor Prof. Gideon Schechtman for motivating me to
prove Proposition \ref{prop:gamma-equivalence}, sharing his
knowledge, and for reading this manuscript.


\section{Duality of Entropy} \label{sec:1}

As emphasized in the Introduction, the proof of Theorem
\ref{thm:main} is immediate once we recall two previously known
results. The first is the recently observed ``Reduction Lemma"
(\cite[Proposition 12]{AMS-Duality-For-Ball}), which uses a
clever iteration procedure to telescopically expand and reduce
the appearing terms. We carefully formulate it below:

\begin{thm}[\cite{AMS-Duality-For-Ball}]
Let $T$ be a convex symmetric body in a Euclidean space such that,
for some constants $a,b \geq 1$, for any convex symmetric body $K
\subset 4T$, one has:
\[
\log N(K, T) \leq b \log N(T^\circ, a^{-1} K^\circ).
\]
Then for any convex symmetric body $K$:
\[
\log N(K, T) \leq b \log_2(48 a) \log N(T^\circ, (8a)^{-1}
K^\circ).
\]
Dually, if $K$ is fixed and the hypothesis holds for all $T$
verifying $K \subset 4T$, then the conclusion holds for any $T$.
\end{thm}

The second known result goes back to the work of \cite{BPST}. It
uses the so called Maurey's Lemma, which (roughly speaking)
estimates the covering number of the convex hull of $m$ points by
the unit-ball of a $K$-convex space. We combine Theorem 6 and the
subsequent remark from \cite{BPST} into the following:

\begin{thm}[\cite{BPST}] \label{thm:BPST}
Let $K,T$ be two symmetric convex bodies in $\Real^n$, such that
$K \subset 4 T$. Then:
\[
\log N(K,T) \leq V \log N(T^\circ,V^{-1} K^\circ),
\]
where $V = \min(V(K),V(T))$ and $V(L)$ is given by
(\ref{eq:defn-V}).
\end{thm}

Combining these two results, we immediately deduce Theorem
\ref{thm:main}. Note that if $V = V(T)$ in Theorem
\ref{thm:main}, we proceed by fixing $T$, applying Theorem
\ref{thm:BPST} for all $K$ satisfying $K \subset 4T$ and use the
first part of the Reduction Lemma to deduce (\ref{eq:final}) for
all $K$; if $V = V(K)$, we fix $K$ and repeat the argument by
interchanging the roles of $K$ and $T$ and using the second part
of the Reduction Lemma.

\medskip

We remark that the proof of Theorem \ref{thm:BPST} in fact gives
an explicit expression for $V$, rather than the implicit one used
in (\ref{eq:defn-V}):
\begin{equation} \label{eq:BPST-V}
V := C_1 \inf \set{\log(C_2 \gamma) (10 T_p(X_B))^q ; K \subset
\gamma B , B \subset 4T, \gamma \geq 1},
\end{equation}
where the infimum runs over all symmetric convex bodies $B$ in
$\Real^n$, $T_p(X_B)$ is the type $p$ ($1<p\leq 2$) constant of
the Banach space $X_B$ whose unit-ball is $B$, $q = p^* =
p/(p-1)$ and $C_1,C_2\geq 1$ are two universal constants (see
\cite{Milman-Schechtman-Book} for the definition of type). Theorem
\ref{thm:main} was formulated using an implicit function $f$ of
$\K(X_B)$, since by several important results of Pisier
(\cite{Pisier-Type-1-Spaces},\cite{Pisier-Type-Implies-K-Convex}),
an infinite dimensional Banach Space is $K$-convex iff it has
some non-trivial type $p>1$. We comment that in
\cite{Pisier-Type-Implies-K-Convex}, an explicit formula bounding
$\K(X_B)$ as a function of $T_p(X_B)$ and $q$ was obtained. It is
possible to obtain an explicit reverse bound using the results in
\cite{Pisier-Stable-Type}, but it is much easier to use an
abstract argument which infers the existence of a $p>1$,
depending solely on $\K(X_B)$, such that $T_p(X_B)$ depends
solely on $\K(X_B)$ (see, e.g. \cite[Lemma
4.2]{KlartagEMilman-2-Convex}). The advantage of using the
$K$-convexity constant $\K(X_B)$ (instead of $T_p(X_B)$ and $q$),
lies in the fact that we may use duality and deduce the other
side of the duality inequality (\ref{eq:other-side}) with the
same $V$, as explained in the Introduction. We also remark that
once $V$ in (\ref{eq:BPST-V}) is expressed using $\K(X_B)$, it is
clear that $V \leq \min(V(K),V(T))$ where $V(L)$ is given by
(\ref{eq:defn-V}). We need this ``separable" estimate on $V$, so
that we may apply the Reduction Lemma (where the estimate on one
of the bodies must be fixed).

\medskip

It is important to note that the proof of Theorem \ref{thm:BPST}
actually connects the notions of covering and convexified packing,
mentioned in the Introduction. For two symmetric convex bodies
$K$ and $T$, the convexified packing number, or convex separation
number, was defined in \cite{AMST-Duality-For-K-Convex} as:
\[
\hat{M}(K,T) = \max \set{N; \begin{array}{c} \exists
x_1,\ldots,x_N \in K \text{ such that } \\ (x_j + \text{int}T)
\cap \text{conv}\set{x_i; i<j} = \emptyset \end{array}}.
\]
Here $\text{int}(T)$ denotes the interior of the set $T$. Note
that we always have $\hat{M}(K,T) \leq N(K,T/2)$ by a standard
argument (see \cite{AMST-Duality-For-K-Convex}). Then the proof
actually shows:

\begin{thm}[\cite{BPST}] \label{thm:main2}
Under the same conditions as in Theorem \ref{thm:BPST}:
\[
\log N(K,T) \leq V \log \hat{M}(K,V^{-1} T).
\]
\end{thm}

Using John's Theorem as in Corollary \ref{cor:main}, we have:

\begin{cor}[\cite{BPST}] \label{cor:main2}
Let $K,T$ be two symmetric convex bodies in $\Real^n$, such that
$K \subset 4 T$. Then:
\[
\log N(K,T) \leq C \log(1+n) \log \hat{M}(K,(C \log(1+n))^{-1} T).
\]
\end{cor}

We mention these variants of Theorem \ref{thm:main} and Corollary
\ref{cor:main} here, because the framework developed in
\cite{AMST-Duality-For-K-Convex} suggests that this is the correct
way to understand the duality problem. The cost of transition
from covering to convex separation, as given by Theorem
\ref{thm:main2} and Corollary \ref{cor:main2}, is a certain
measure of the complexity of the bodies involved. Once the
transition is achieved, the duality framework developed in
\cite{AMS-Duality-For-Ball} and \cite{AMST-Duality-For-K-Convex}
finishes the job. Indeed, it was shown in
\cite{AMST-Duality-For-K-Convex} that the convex separation
numbers always satisfy a duality relation, for any pair of
symmetric convex bodies $K,T$:
\[
\hat{M}(K,T) \leq \hat{M}(T^\circ,K^\circ/2)^2.
\]
Using Theorem \ref{thm:main2}, we conclude that when $K \subset
4T$:
\begin{eqnarray}
\nonumber \log N(K,T) &\leq& V \log \hat{M}(K,V^{-1} T) \leq  2 V
\log \hat{M}(T^\circ,(2V)^{-1} K^\circ) \\
\nonumber &\leq&  2 V \log N(T^\circ,(4V)^{-1} K^\circ).
\end{eqnarray}
The Reduction Lemma now immediately gives Theorem \ref{thm:main}.

\medskip

To conclude this section, we mention that Theorem \ref{thm:main2}
is already stronger than all of the results in
\cite{AMST-Duality-For-K-Convex} connecting the covering and the
convex separation numbers. The technique involving the use of
Maurey's Lemma, which was also used in
\cite{AMST-Duality-For-K-Convex} (see also \cite{Shiri-PhD}), is
optimally exploited in the proof of Theorem \ref{thm:main2}
(Theorem 6 in \cite{BPST}), by using a clever iteration procedure,
producing the $\log$ factor in the various definitions
(\ref{eq:defn-V}) and (\ref{eq:BPST-V}) of $V$. All previous
general results (with no restriction on $K$ and $T$) pay a linear
penalty in the Banach-Mazur distance from ``low-complexity"
bodies, which may be as large as $\sqrt{n}$.


\section{Duality of Talagrand's $\gamma_p$ Functionals}
\label{sec:2}

Given Corollary \ref{cor:main}, proving Theorem \ref{thm:gamma}
is rather elementary, although it involves an analogue to
Sudakov's Minoration bound which we have not encountered
elsewhere. Before proceeding, we remark that for our purposes, it
is totally immaterial whether the points $\set{x_i}$ in the
definition of $N(K,T)$ are chosen to lie inside $K$ or not.
Indeed, denoting by $N'(K,T)$ the variant where the points
\emph{are} required to lie inside $K$, it is elementary to check
that:
\[
N'(K,2T) \leq N(K,T) \leq N'(K,T).
\]
Since throughout this note we allow the insertion of homothety
constants in all expressions, or multiplying the entropy numbers
by universal constants, this lack of distinction is well
justified.

\medskip

First, recall that by Dudley's entropy bound
(\cite{Talagrand-Book}):
\begin{equation} \label{eq:Dudley}
\gamma_p(K,T) \leq C_p \sum_{k \geq 1} k^{1/p-1} e_k(K,T),
\end{equation}
where $C_p>0$ is some constant depending on $p$. The argument is
elementary:
\[
\gamma_p(K,T) := \inf \sup_{x \in K} \sum_{j \geq 0} 2^{j/p}
d_T(x,K_j) \leq \inf \sum_{j \geq 0} 2^{j/p} \sup_{x \in K}
d_T(x,K_j).
\]
Choosing $K_j$ to be the set of $2^{2^j}$ points (inside $K$)
attaining the minimum in the definition of $N(K,e_{2^j}(K,T) T)$,
we see that:
\[
\gamma_p(K,T) \leq \sum_{j \geq 0} 2^{j/p} e_{2^j}(K,T).
\]
It is elementary to verify that for $p\geq 1$ and $j\geq 1$:
\[
2^{j/p} \leq C_p \sum_{k=2^{j-1}}^{2^j-1} k^{1/p-1},
\]
where $C_p = (p(1 - 2^{-1/p}))^{-1}$. Since $e_k$ is a
non-increasing sequence, we have:
\begin{eqnarray}
\nonumber \gamma_p(K,T) \leq e_1(K,T) + C_p \sum_{j \geq 1}
\sum_{k=2^{j-1}}^{2^j-1} k^{1/p-1} e_{2^j}(K,T) \\
\nonumber \leq e_1(K,T) + C_p \sum_{k \geq 1} k^{1/p-1} e_k(K,T)
\leq (1+C_p) \sum_{k \geq 1} k^{1/p-1} e_k(K,T).
\end{eqnarray}
A similar argument works for $0<p<1$.

Dudley's entropy upper bound appears naturally when studying the
supremum of Gaussian processes on a set $K$, e.g. $E \sup_{x \in
K} \scalar{x,G}$ where $G$ is a Gaussian r.v. As mentioned in the
Introduction, a deep theorem of Talagrand asserts that the latter
expectancy is in fact equivalent (to within universal constants)
to $\gamma_2(K,D)$ where $D$ is an ellipsoid corresponding to the
covariance of $G$. The corresponding lower bound on $E \sup_{x
\in K} \scalar{x,G}$ is due to Sudakov (\cite{Sudakov}):
\[
E \sup_{x \in K} \scalar{x,G} \geq c \sup_{k \geq 1} k^{1/2}
e_k(K,D).
\]
When the body $T$ is not an ellipsoid or when $p\neq 2$, there is
no direct connection between Gaussian processes and
$\gamma_p(K,T)$. Nevertheless, we note that the analogue to
Sudakov's Minoration bound holds in full generality:

\begin{lem} \label{lem:Sudakov}
\[
\gamma_p(K,T) \geq 2^{-1/p} \sup_{k \geq 1} k^{1/p} e_k(K,T).
\]
\end{lem}

\begin{proof}
Let $k\geq 1$ be given, and let $j \geq 0$ be such that $2^j \leq
k < 2^{j+1}$. Then:
\[
\gamma_p(K,T) := \inf \sup_{x \in K} \sum_{l \geq 0} 2^{l/p}
d_T(x,K_l) \geq \inf \sup_{x \in K} 2^{j/p} d_T(x,K_j).
\]
Since for any admissible set $K_j$ we have $|K_j| = 2^{2^j} \leq
2^k$, it follows by definition that $\sup_{x \in K} d_T(x,K_j)
\geq e_k(K,T)$. Hence:
\[
\gamma_p(K,T) \geq 2^{j/p} e_k(K,T) \geq 2^{-1/p} k^{1/p}
e_k(K,T).
\]
Since $k \geq 1$ was arbitrary, the assertion follows.
\end{proof}

We will need one last lemma for the proof of Theorem
\ref{thm:gamma}:

\begin{lem} \label{lem:entropy-tail}
For all $k \geq 3n$:
\[
e_k(K,T) \leq 2 e_n(K,T) \exp(-c k/n),
\]
where $c>0$ is some universal constant.
\end{lem}

\begin{proof}
Denote $e_k = e_k(K,T)$ and $e_n = e_n(K,T)$ for short. W.l.o.g.
we assume that $N(K,e_k T) = 2^k$ and $N(K, e_n T) = 2^n$.
Obviously we have:
\begin{equation} \label{eq:triangle}
N(K,e_k T) \leq  N(K, e_n T) N(e_n T , e_k T).
\end{equation}
Also $N(e_n T , e_k T) = N(T, \frac{e_k}{e_n} T)$, and by a
standard volume estimation argument, we can find an $e_k/e_n$
$T$-net of $T$ with cardinality no greater than
$(1+\frac{2e_n}{e_k})^n$. Plugging everything into
(\ref{eq:triangle}), we see that:
\[
2^k \leq 2^n \brac{1+\frac{2e_n}{e_k}}^n,
\]
or equivalently:
\[
\exp\brac{ \log(2) \frac{k-n}{n}} - 1 \leq \frac{2 e_n}{e_k}.
\]
Since $k \geq 3n$, we can find a universal constant $c>0$ such
that:
\[
\exp\brac{\log(2) \frac{k-n}{n}} - 1 \geq \exp\brac{c
\frac{k}{n}}.
\]
The assertion now readily follows.
\end{proof}

We can now deduce the following equivalence, up to a $log(n)$
term, of the $\gamma_p$ functional, Sudakov's lower bound and
Dudley's upper bound. Although this is probably known, we did not
find a reference for it, so we include a proof for completeness.

\begin{prop} \label{prop:gamma-equivalence}
Let $K,T$ denote two symmetric convex bodies in $\Real^n$, and
denote $e_k = e_k(K,T)$ and $\gamma_p = \gamma_p(K,T)$ for short.
Then for any $p>0$:
\[
2^{-1/p} \sup_{k \geq 1} k^{1/p} e_k \leq \gamma_p \leq C_p
\sum_{k \geq 1} k^{1/p-1} e_k \leq \log(1+n) C'_p \sup_{k \geq 1}
k^{1/p} e_k,
\]
where $C_p,C'_p > 0$ are universal constants depending solely on
$p$.
\end{prop}

\begin{proof}
The first inequality is Sudakov's lower bound (Lemma
\ref{lem:Sudakov}) and the second one is Dudley's upper bound
(\ref{eq:Dudley}). We will show the third inequality. Let us
split the sum $\sum_{k \geq 1} k^{1/p-1} e_k$ into two parts, up
to and from $k=3n$. For the first part, we obviously have:
\[
\sum_{k=1}^{3n-1} k^{1/p-1} e_k \leq \sum_{k=1}^{3n-1}
\frac{1}{k} \; \sup_{k\geq 1} k^{1/p} e_k \leq C \log(1+n)
\sup_{k\geq 1} k^{1/p} e_k.
\]
We use Lemma \ref{lem:entropy-tail} to evaluate the second sum:
\[
\sum_{k\geq 3n} k^{1/p-1} e_k \leq 2 e_n \sum_{k\geq 3n}
k^{1/p-1} \exp(-c k/n).
\]
For $p\geq 1$, $k^{1/p-1}$ is non-increasing, so we use:
\begin{eqnarray}
\nonumber \sum_{k\geq 3n} k^{1/p-1} \exp(-c k/n) \leq (3n)^{1/p-1}
\sum_{k\geq 3n} \exp(-c k/n) \\
\nonumber = (3n)^{1/p-1} \frac{\exp(-3c)}{1-\exp(-c/n)} \leq
(3n)^{1/p-1} \exp(-3c) \frac{n}{c'} \leq C n^{1/p}.
\end{eqnarray}
For $0<p<1$, we evaluate the sum with an integral (although the
series may not be monotone, is has at most one extremal point,
and this can be handled by a loose estimate):
\begin{eqnarray}
\nonumber \sum_{k\geq 3n} k^{1/p-1} \exp(-c k/n) \leq 3
\int_{3n-1}^\infty x^{1/p-1} \exp(-c x/n) dx \\
\nonumber \leq 3 \brac{\frac{n}{c}}^{1/p} \int_0^\infty x^{1/p-1}
\exp(-x) dx = 3 c^{-1/p} \Gamma(1/p) n^{1/p}.
\end{eqnarray}
We conclude that in both cases:
\[
\sum_{k\geq 3n} k^{1/p-1} e_k \leq C'_p n^{1/p} e_n \leq C'_p
\sup_{k \geq 1} k^{1/p} e_k.
\]
Summing the two parts together, we conclude the proof.
\end{proof}

Using Corollary \ref{cor:main}, the proof of Theorem
\ref{thm:gamma} is now clear:

\begin{proof}[Proof of Theorem \ref{thm:gamma}]
Corollary \ref{cor:main} implies that:
\[
e_k(K,T) \leq C \log(1+n) e_{k / (C \log(1+n) \log \log(2+n))}
(T^\circ,K^\circ),
\]
for some universal constant $C \geq 1$ and all $k\geq 0$. Using
Proposition \ref{prop:gamma-equivalence} twice, we conclude:
\begin{eqnarray}
\nonumber \gamma_p(K,T) &\leq& C'_p \log(1+n) \sup_{k \geq 1}
k^{1/p} e_k(K,T) \\
\nonumber &\leq& C''_p \log(1+n)^2 (\log(1+n) \log
\log(2+n))^{1/p} \sup_{k \geq 1} k^{1/p} e_k(T^\circ,K^\circ) \\
\nonumber &\leq& C_p \log(1+n)^{2+1/p} \log \log(2+n)^{1/p}
\gamma_p(T^\circ,K^\circ)
\end{eqnarray}

\end{proof}

\bibliographystyle{amsalpha}
\bibliography{../../ConvexBib}

\end{document}